\documentclass[11pt]{article}
\usepackage{latexsym,amsmath,amsfonts}
\DeclareSymbolFont{lettersA}{U}{pxmia}{m}{it}
\DeclareMathSymbol{\piup}{\mathord}{lettersA}{"19}
\begin{document}
\begin{center}
{\LARGE\textbf{Laws of Large Numbers of Subgraphs in\\
Directed Random Geometric Networks}}\\
\bigskip
Yilun Shang\footnote{Department of Mathematics, Shanghai Jiao Tong
University, 800
DongChuan  Road, Shanghai, CHINA. E-mail address: shyl@sjtu.edu.cn}\\
\medskip
\end{center}

\begin{abstract}
Given independent random points $\mathcal{X}_n=\{X_1,\cdots,X_n\}$
in $\mathbb{R}^2$, drawn according to some probability density
function $f$ on $\mathbb{R}^2$, and a cutoff $r_n>0$ we construct a
random geometric digraph $G(\mathcal{X}_n,\mathcal{Y}_n,r_n)$ with
vertex set $\mathcal{X}_n$. Each vertex $X_i$ is assigned uniformly
at random a sector $S_i$, of central angle $\alpha$ with inclination
$Y_i$, in a circle of radius $r_n$ (with vertex $X_i$ as the
origin). An arc is present from $X_i$ to $X_j$, if $X_j$ falls in
$S_i$. We also introduce another random geometric digraph
$G(\mathcal{X}_n,\mathcal{R}_n)$ with vertex set
$\mathcal{X}_n=\{X_1,\cdots,X_n\}$ in $\mathbb{R}^d$, $d\ge1$ and an
arc present from  $X_i$ to $X_j$ if $||X_i-X_j||<R_{n,i}$. Here
$\{R_{n,i}\}_{i\ge1}$ are i.i.d. random variables and we may take an
arbitrary norm $||\cdot||$. In this paper we investigate two kinds
of small subgraphs---induced and isolated---in the above two
directed networks, which contribute to understanding the local
topology of many spatial networks, such as wireless communication
networks. We give some strong laws of large numbers of subgraph
counts thus extending those results of Penrose [Random Geometric
Graphs, Oxford University Press, 2003].


\end{abstract}

\bigskip
\normalsize

\noindent{\Large\textbf{1. Introduction}}

\smallskip
In the last decade there has been a resurgence of interest in the
analysis of random geometric graphs (RGGs) particularly in the
context of ad hoc wireless networks. An elegant written tutorial of
random geometric graph theory is available in \cite{2}, and the
paper \cite{1} is a more recent survey emphasizing wireless
networks. An RGG is usually constructed as follows. Let $||\cdot||$
be some norm on $\mathbb{R}^d$, and $r_n$ be a real sequence. Let
$\mathcal{X}_n=\{X_1, X_2,\cdots,X_n\}$, $\{X_i\}$ are i.i.d. random
$d$-vectors in $\mathbb{R}^d$ having a common probability density
function $f$. We denote by $G(\mathcal{X}_n,r_n)$ the graph with
vertex set $\mathcal{X}_n$ and with an edge $X_iX_j$ if and only if
$||X_i-X_j||\le r_n$ for $i\not=j$. Note that $G(\mathcal{X}_n,r_n)$
is isotropic and thus undirected, which is less appropriate in many
practical applications such as wireless sensor networks. The issue
of small subgraph counts are dealt with here for two kinds of
directed models of RGG (for formal definitions see below). In the
case of wireless networks, we are sometimes interested to know the
number of a desirable local configuration involving a small number
of transmitters and receivers. This is especially true if there are
a small number of nodes with special capabilities, e.g. data
collection centers in sensor network process the data collected by
the beacon nodes that help in self-organization of the network. The
small subgraph counts are also of independent interest in random
graph theory in various guises. The concerned results on small
subgraph for classical Erd\"os-R\'enyi random graphs are discussed
in detail in \cite{10}(chap. 4) and \cite{11}(chap. 3), and for
asymptotic results in random geometric graphs, see \cite{2}(chap. 2)
and \cite{50}, while for exact formulae treated in the circumstance
of wireless network, see \cite{7}.

In this short paper, we extend the method of Penrose \cite{2} and
establish some strong laws of large numbers of small subgraphs in
directed geometric networks in some limiting regimes. We now define
two random geometric digraphs to be used in this work. The archetype
(with uniformly distributed points in $[0,1]^2$) of the first one
has been proposed in \cite{3}, called ``random scaled sector
graph'', to model the ``Small Dust'' sensor networks using optical
communication. Some graph theoretic properties have been addressed
for this model, see e.g. \cite{3,4,8,6,20}, mainly using
combinatorial techniques.

\smallskip
\noindent\textbf{Definition 1.} \itshape Let $||\cdot||$ be
Euclidean norm equipped on $\mathbb{R}^2$. Let $\alpha\in
(0,2\piup]$. Let $\mathcal{X}_n=\{X_1, X_2,\cdots,X_n\}$ be i.i.d.
random vectors with a common density function
$f:\mathbb{R}^2\rightarrow\mathbb{R}$. Let
$\mathcal{Y}_n=\{Y_1,Y_2,\cdots,Y_n\}$ be i.i.d. random variables,
uniformly distributed on $[0,2\piup)$. Associate every point $X_i\in
\mathcal{X}_n$ a sector, which is centered at $X_i$, with radius
$r_n$, amplitude $\alpha$ and elevation $Y_i$ with respect to the
horizontal direction anticlockwise. This sector is denoted as
$S(X_i,Y_i,r_n)$. We denote by $G(\mathcal{X}_n,\mathcal{Y}_n,r_n)$
the digraph with vertex set $\mathcal{X}_n$, and with an arc
$(X_i,X_j)$, $i\not=j$, presents if and only if $X_j\in
S(X_i,Y_i,r_n)$. \normalfont

\smallskip
For technical reasons we always assume $r_n\rightarrow0$ as
$n\rightarrow\infty$. We also assume that $f$ is bounded and a.s.
continuous throughout the paper. We mention that the above
assumptions imposed on $f$ is rather mild; in fact typical
distributions such as normal distribution and $f=1_{[0,1]^d}$ are
clearly satisfied. Now we introduce another model
$G(\mathcal{X}_n,\mathcal{R}_n)$ with random cutoffs motivated by
Boolean model in continuum percolation \cite{12}. For the sake of
convenience, we still choose to employ the signs $\mathcal{X}_n$,
$f$ and $||\cdot||$ with a little ambiguous (see below), however,
the right meaning will be clear in the context and no confusion will
be incurred.

\smallskip
\noindent\textbf{Definition 2.} \itshape Let $||\cdot||$ be any norm
equipped on $\mathbb{R}^d$, $d\ge 1$. Let $\mathcal{X}_n=\{X_1,
X_2,\cdots,X_n\}$ be i.i.d. random vectors with a common density
function $f:\mathbb{R}^d\rightarrow\mathbb{R}$. Let $R_n$ be a
positive random variable with probability distribution function
$Q_n$ and density function $q_n$. For each point $X_i$, we associate
a ball $B(X_i,R_{n,i})$ with radius $R_{n,i}$, centered at $X_i$,
independent of other points. $\{R_{n,i}\}_{i=1}^{n}$ are independent
copies of $R_n$ and set $\mathcal{R}_n:=\{R_{n,1}, R_{n,2},\cdots,
R_{n,n}\}$. We denote by $G(\mathcal{X}_n,\mathcal{R}_n)$ the
digraph with vertex set $\mathcal{X}_n$, and with an arc $(X_i,X_j)$
originating from $X_i$ and terminating in $X_j$ if and only if
$X_j\in B(X_i,R_{n,i})$.\normalfont

\smallskip
As usual, we shall impose a certain decaying condition on $R_n$.
Here we assume
$ER_n^d=\int_0^{\infty}r^d\mathrm{d}Q_n(r)\rightarrow0$ as
$n\rightarrow\infty$, throughout the paper. We will investigate two
kinds of subgraphs in the above two models; one is induced subgraph
and the other is isolated subgraph. Induced subgraph is defined in
its usual meaning, see e.g.\cite{13}. Suppose $G$ is a digraph, a
subgraph $H$ is isolated in $G$ if $H$ is an induced subgraph and
there are no arcs leave $H$. If $G$ is undirected, then a connected
isolated subgraph of $G$ is just a component.

Before going further we will need some other definitions. Given a
finite set $\mathcal{X}$ in $\mathbb{R}^d$, let card($\mathcal{X}$)
denote the number of points in $\mathcal{X}$ and let $|\cdot|$ be
the $d$-dimensional Lebesgue measure, which is easy to discriminate
from the sign for absolute value in the context. In the rest of the
paper, let $f_{\max}:=\sup\{t:|\{f(x)>t\}|>0\}$ be the essential
supremum of the probability density function $f$. As mentioned
before, we assume $f_{\max}<\infty$. For a set $A\in\mathbb{R}^d$,
let $\mathcal{X}(A)$ denote the number of points of $\mathcal{X}$
located in $A$. Denote $D(0,1)$ as the unit disk in $\mathbb{R}^2$,
then $|D(0,1)|=\piup$ and also set $\theta:=|B(0,1)|$ w.r.t some
given norm. Let $C$, $C'$ etc. be various positive constants, and
the values may change from line to line.

The rest of this paper is organized as follows. Section 2 contains
the statement of main results and proofs are provided in Section 3.
We finally draw conclusions in Section 4.

\bigskip
\noindent{\Large\textbf{2. Statement of main results}}

\smallskip
For $k\in\mathbb{N}$, let $T$ be a fixed connected graph on $k$
vertices. We say $T$ is feasible if either
$P(G(\mathcal{X}_k,\mathcal{Y}_k,r)\cong T)>0$ for some $r>0$ or
$P(G(\mathcal{X}_k,\mathcal{R}_n)\cong T)>0$ for some $\{r_{n,1},
r_{n,2},\cdots,r_{n,k}\}$. Let $H_{n,T}$ and $\tilde{H}_{n,T}$ be
the number of induced subgraphs and isolated subgraphs of
$G(\mathcal{X}_n,\mathcal{Y}_n,r_n)$ isomorphic to $T$
($T$-subgraphs) respectively. Likewise, let $G_{n,T}$ and
$\tilde{G}_{n,T}$ be the number of induced and isolated
$T$-subgraphs of $G(\mathcal{X}_n,\mathcal{R}_n)$ respectively. For
a finite set $\mathcal{X}\subseteq\mathbb{R}^2$ and a point
$\mathcal{Y}\in[0,2\piup)^{\mathrm{card}(\mathcal{X})}$, we define
indicator random variables
$h_T(\mathcal{X},\mathcal{Y}):=1_{[G(\mathcal{X},\mathcal{Y},1)\cong
T]}$ and
$h_{n,T}(\mathcal{X},\mathcal{Y}):=1_{[G(\mathcal{X},\mathcal{Y},r_n)\cong
T]}$. For a finite set $\mathcal{X}\subseteq\mathbb{R}^d$ and
$\mathcal{R}_n=\{R_{n,1}, R_{n,2},\cdots,
R_{n,\mathrm{card}(\mathcal{X})}\}$, we define
$g_{n,T}(\mathcal{X},\mathcal{R}_n):=1_{[G(\mathcal{X},\mathcal{R}_n)\cong
T]}$.

The basic tool we shall need in the proofs is the following Azuma's
inequality; we refer the reader to \cite{9} for a proof and a wealth
of materials regarding that topic.

\smallskip
\noindent\textbf{Lemma 1.}\quad\itshape Suppose $M_1,
M_2,\cdots,M_n$ is a martingale with corresponding martingale
difference sequence $D_1, D_2,\cdots,D_n$, where $D_i:=M_i-M_{i-1}$
$(2\le i\le n)$ and $D_1:=M_1-EM_1$. Then for any $a>0$, we have
$$
P\Big(\Big|\sum_{i=1}^nD_i\Big|>a\Big)\le
2\exp\Big(-\frac{a^2}{2\sum_{i=1}^n||D_i||^2_{\infty}}\Big),
$$
where $||D_i||^2_{\infty}:=\inf\{b:P(|D_i|\le b)=1\}$. \normalfont

\smallskip
In the sequel we sometimes use a generalized version with
``tolerance'' of Azuma's inequality \cite{5}:

\smallskip
\noindent\textbf{Lemma 2.}\itshape (\cite{5})\quad Suppose $M_1,
M_2,\cdots,M_n$ is a martingale with corresponding martingale
difference sequence $D_1, D_2,\cdots,D_n$, where $D_i:=M_i-M_{i-1}$
$(2\le i\le n)$ and $D_1:=M_1-EM_1$. Then for any $a,b>0$,
$$
P\Big(\Big|\sum_{i=1}^nD_i\Big|>a\Big)\le
2\exp\Big(-\frac{a^2}{32nb^2}\Big)+\Big(1+\frac{2\sup_i||D_i||_{\infty}}{a}\Big)\sum_{i=1}^nP\big(|D_i|>b\big).
$$
\normalfont

\smallskip
By the notations defined in the beginning of this section, we have
the following strong laws of large numbers under various regimes:
\smallskip

\noindent\textbf{Theorem 1.}\quad\itshape Suppose $T$ is a connected
feasible graph of order $k$, $k\in\mathbb{N}$. Let
$nr_n^2\rightarrow\lambda\in(0,\infty)$. Then
$$
\lim_{n\rightarrow\infty}n^{-1}\tilde{H}_{n,T}=k^{-1}\int_{\mathbb{R}^2}\varphi_T(\lambda
f(x))f(x)\mathrm{d}x\quad a.s.
$$
where
$$\varphi_T(t):=\left\{
\begin{array}{lr}
\frac{t^{k-1}}{(k-1)!(2\piup)^k}\underbrace{\int_0^{2\piup}\cdots\int_0^{2\piup}}_{k}\underbrace{\int_{\mathbb{R}^2}\cdots\int_{\mathbb{R}^2}}_{k-1}h_T(0,x_2,\cdots,x_k,y_1,\cdots,y_k)\\
\qquad\qquad\cdot
e^{-t|S(0,x_2,\cdots,x_k,y_1,\cdots,y_k)|}\mathrm{d}x_2\cdots\mathrm{d}x_k\mathrm{d}y_1\cdots\mathrm{d}y_k&
k\ge2;\\
e^{-\frac{t\alpha}{2}} & k=1,
\end{array}
\right.
$$
and $S(x_1,\cdots,x_k,y_1,\cdots,y_k):=\cup_{i=1}^kS(x_i,y_i,1)$.
\normalfont

\medskip
\noindent\textbf{Theorem 2.}\quad\itshape Let $T$ be a single
vertex. Suppose $nR_n^d\rightarrow\lambda\in[0,\infty)$ in
probability. Then
$$
\lim_{n\rightarrow\infty}n^{-1}\tilde{G}_{n,T}=\int_{\mathbb{R}^d}e^{-\theta\lambda
f(x)}f(x)\mathrm{d}x\quad a.s.
$$
\normalfont

\medskip
\noindent\textbf{Theorem 3.}\quad\itshape For $k\ge2$, let $T$ be a
connected feasible graph of order $k$. Suppose $nr_n^2\rightarrow0$
and $\ln n=o\big(n^{2k-1}r_n^{4(k-1)}\big)$, as
$n\rightarrow\infty$. Then
$$
\lim_{n\rightarrow\infty}n^{-k}r_n^{-2(k-1)}\tilde{H}_{n,T}=\mu_T\quad
a.s.
$$
where
\begin{multline}
\mu_T:=\frac1{k!(2\piup)^k}\underbrace{\int_0^{2\piup}\cdots\int_0^{2\piup}}_{k}\underbrace{\int_{\mathbb{R}^2}\cdots\int_{\mathbb{R}^2}}_{k}h_T(0,z_2,\cdots,z_k,y_1,\cdots,y_k)\\
\cdot
f^k(x)\mathrm{d}x\mathrm{d}z_2\cdots\mathrm{d}z_k\mathrm{d}y_1\cdots\mathrm{d}y_k.\nonumber
\end{multline}
\normalfont

\medskip
\noindent\textbf{Theorem 4.}\quad\itshape For $k\ge2$, let $T$ be a
connected feasible graph of order $k$. Suppose $f$ has bounded
support set (denoted by $suppf$). Suppose $r_n\rightarrow0$ and
there is a constant $\delta>0$ such that $\liminf
n^{2k-1-\delta}r_n^{4(k-1)}>0$. Then
$$
\lim_{n\rightarrow\infty}n^{-k}r_n^{-2(k-1)}H_{n,T}=\mu_T \quad a.s.
$$
\normalfont

\bigskip
\noindent{\Large\textbf{3. Proofs}}

\smallskip
To prove our main theorems, we first derive several asymptotic
results (Propositions 1$\sim$4) for the means of subgraph counts
$H_{n,T}$, $\tilde{H}_{n,T}$ and $\tilde{G}_{n,T}$.

\medskip
\noindent\textbf{Proposition 1.}\quad\itshape For $k\ge2$, let $T$
be a connected feasible graph of order $k$. Suppose
$r_n\rightarrow0$ as $n\rightarrow\infty$. Then
$$
\lim_{n\rightarrow\infty}n^{-k}r_n^{-2(k-1)}EH_{n,T}=\mu_T
$$
where $\mu_T$ is defined in Section 2.
\normalfont

\medskip
\noindent\textbf{Proof}. From the linearity of expectation,
$EH_{n,T}={n \choose k}Eh_{n,T}(\mathcal{X}_k,\mathcal{Y}_k)$.
Thereby
\begin{eqnarray}
EH_{n,T}&\hspace{-5pt}=&\hspace{-5pt}\frac1{(2\piup)^k}{n \choose
k}\int_0^{2\piup}\cdots\int_0^{2\piup}\int_{\mathbb{R}^2}\cdots\int_{\mathbb{R}^2}h_{n,T}(x_1,\cdots,x_k,y_1,\cdots,y_k)\nonumber\\
 & &\cdot f^k(x_1)\mathrm{d}x_1\cdots\mathrm{d}x_k\mathrm{d}y_1\cdots\mathrm{d}y_k\nonumber\\
 & &\hspace{-5pt}+\frac1{(2\piup)^k}{n \choose
k}\int_0^{2\piup}\cdots\int_0^{2\piup}\int_{\mathbb{R}^2}\cdots\int_{\mathbb{R}^2}h_{n,T}(x_1,\cdots,x_k,y_1,\cdots,y_k)\nonumber\\
 & &\cdot\bigg(\prod_{i=1}^kf(x_i)-f^k(x_1)\bigg)\mathrm{d}x_1\cdots\mathrm{d}x_k\mathrm{d}y_1\cdots\mathrm{d}y_k\label{1}
\end{eqnarray}
Let $x_1=x$ and $x_i=x_1+r_nz_i$ for $2\le i\le k$, then the first
term on the right hand side of (\ref{1}) equals
\begin{multline}
\frac1{(2\piup)^k}{n \choose
k}r_n^{2(k-1)}\int_0^{2\piup}\cdots\int_0^{2\piup}\int_{\mathbb{R}^2}\cdots\int_{\mathbb{R}^2}h_{n,T}(x,x+r_nz_2,\cdots,x+r_nz_k,y_1,\cdots,y_k)\\
\cdot
f^k(x)\mathrm{d}x\mathrm{d}z_2\cdots\mathrm{d}z_k\mathrm{d}y_1\cdots\mathrm{d}y_k.\nonumber
\end{multline}
We have
$h_{n,T}(x,x+r_nz_2,\cdots,x+r_nz_k,y_1,\cdots,y_k)=h_T(0,z_2,\cdots,z_k,y_1,\cdots,y_k)$.
Since $f$ is bounded and $T$ is a connected graph, $f^k(x)$ and
$h_T(0,z_2,\cdots,z_k,y_1,\cdots,y_k)$ are integrable on
$\mathbb{R}^2$ and $(\mathbb{R}^2)^{k-1}\times[0,2\piup)^k$
respectively. Hence, the first term on the right hand side of
(\ref{1}) tends to $n^{k}r_n^{2(k-1)}\mu_T$ by the dominated
convergence theorem.

Set
$\eta_n(x_1):=\int_{D(x_1,kr_n)}\cdots\int_{D(x_1,kr_n)}r_n^{-2(k-1)}\big|\prod_{i=2}^kf(x_i)-f^{k-1}(x_1)\big|\mathrm{d}x_2\cdots\mathrm{d}x_k$.
It's easy to see that the absolute value of the second term on the
right hand side of (\ref{1}) multiplied by $n^{-k}r_n^{-2(k-1)}$ is
bounded by $\int_{\mathbb{R}^2}\eta_n(x_1)f(x_1)\mathrm{d}x_1$. If
$x_1$ is a continuous point of $f$, then $\eta_n(x_1)\rightarrow0$
as $n\rightarrow\infty$ by the mean value theorem for integrals.
Therefore, by the dominated convergence theorem and the assumption
of almost everywhere continuity of $f$, we have
$\int_{\mathbb{R}^2}\eta_n(x_1)f(x_1)\mathrm{d}x_1\rightarrow0$.
This concludes the proof. $\Box$

\medskip
\noindent\textbf{Proposition 2.}\quad\itshape For $k\ge2$, let $T$
be a connected feasible graph of order $k$. Suppose
$nr_n^2\rightarrow0$ as $n\rightarrow\infty$. Then
$$
\lim_{n\rightarrow\infty}n^{-k}r_n^{-2(k-1)}E\tilde{H}_{n,T}=\mu_T
$$
where $\mu_T$ is defined in Section 2. \normalfont

\medskip
\noindent\textbf{Proof}. Let
$\tilde{h}_{n,T}(\mathcal{X},\mathcal{Y}):=1_{[G(\mathcal{X},\mathcal{Y},r_n)\
\mathrm{is}\ \mathrm{isolated}\ \mathrm{and}\ \cong T]}$. Therefore,
\begin{eqnarray*}
E\tilde{H}_{n,T}&=&{n \choose
k}E\tilde{h}_{n,T}(\mathcal{X}_k,\mathcal{Y}_k)\\
 &=&{n \choose k}Eh_{n,T}(\mathcal{X}_k,\mathcal{Y}_k)\cdot P(G(\mathcal{X}_k,\mathcal{Y}_k,r_n)\
\mathrm{is}\
\mathrm{isolated}|G(\mathcal{X}_k,\mathcal{Y}_k,r_n)\cong T)\\
 &:=&EH_{n,T}\cdot P_1
\end{eqnarray*}
Since $T$ is a connected graph and by the assumed asymptotic
behavior of $r_n$, we get as $n\rightarrow\infty$,
$$
1\ge P_1\ge (1-P(X_1\in D(0,kr_n)))^{n-k}\ge
(1-f_{\max}\piup(kr_n)^2)^{n-k}\rightarrow1.
$$
Consequently, $E\tilde{H}_{n,T}=(1+o(1))EH_{n,T}$. By using
Proposition 1, we complete the proof. $\Box$

\medskip
\noindent\textbf{Proposition 3.}\quad\itshape For $k\in\mathbb{N}$,
let $T$ be a connected feasible graph of order $k$. Suppose
$nr_n^2\rightarrow\lambda\in(0,\infty)$ as $n\rightarrow\infty$.
Then
$$
\lim_{n\rightarrow\infty}n^{-1}E\tilde{H}_{n,T}=k^{-1}\int_{\mathbb{R}^2}\varphi_T(\lambda
f(x))f(x)\mathrm{d}x
$$
where $\varphi_T(\cdot)$ is defined in Section 2. \normalfont

\medskip
\noindent\textbf{Proof}. By the definition of $h_{n,T}$ and similar
with the beginning of the proof of Proposition 1, we have
\begin{eqnarray}
n^{-1}E\tilde{H}_{n,T}&\hspace{-5pt}=&\hspace{-5pt}\frac1{n(2\piup)^k}{n
\choose
k}\int_0^{2\piup}\cdots\int_0^{2\piup}\int_{\mathbb{R}^2}\cdots\int_{\mathbb{R}^2}h_{n,T}(x_1,\cdots,x_k,y_1,\cdots,y_k)\nonumber\\
 & &\cdot\Big(1-\int_{\cup_{i=1}^kS(x_i,y_i,r_n)}f(x)\mathrm{d}x\Big)^{n-k}f^k(x_1)\mathrm{d}x_1\cdots\mathrm{d}x_k\mathrm{d}y_1\cdots\mathrm{d}y_k\nonumber\\
 & &\hspace{-5pt}+\frac1{n(2\piup)^k}{n \choose
k}\int_0^{2\piup}\cdots\int_0^{2\piup}\int_{\mathbb{R}^2}\cdots\int_{\mathbb{R}^2}h_{n,T}(x_1,\cdots,x_k,y_1,\cdots,y_k)\nonumber\\
 & &\cdot\Big(1-\int_{\cup_{i=1}^kS(x_i,y_i,r_n)}f(x)\mathrm{d}x\Big)^{n-k}\nonumber\\
 & &\cdot\bigg(\prod_{i=1}^kf(x_i)-f^k(x_1)\bigg)
\mathrm{d}x_1\cdots\mathrm{d}x_k\mathrm{d}y_1\cdots\mathrm{d}y_k\label{2}
\end{eqnarray}
Let $x_i=x_1+r_nz_i$ for $2\le i\le k$, then the first term on the
right hand side of (\ref{2}) tends to
\begin{multline}
\frac1{n(2\piup)^k}{n \choose
k}r_n^{2(k-1)}\int_0^{2\piup}\cdots\int_0^{2\piup}\int_{\mathbb{R}^2}\cdots\int_{\mathbb{R}^2}h_{n,T}(x_1,x_1+r_nz_2,\cdots,x_1+r_nz_k,y_1,\cdots,y_k)\\
\cdot e^{(n-k)\ln
\big(1-\int_{S(x_1,y_1,r_n)\cup\cup_{i=2}^kS(x_1+r_nz_i,y_i,r_n)}f(x)\mathrm{d}x\big)}
f^k(x_1)\mathrm{d}x_1\mathrm{d}z_2\cdots\mathrm{d}z_k\mathrm{d}y_1\cdots\mathrm{d}y_k\nonumber
\end{multline}
which is further asymptotic to
\begin{multline}
\frac{\lambda^{k-1}}{k!(2\piup)^k}\int_0^{2\piup}\cdots\int_0^{2\piup}\int_{\mathbb{R}^2}\cdots\int_{\mathbb{R}^2}h_T(0,z_2,\cdots,z_k,y_1,\cdots,y_k)\\
\cdot e^{-\lambda
f(x_1)|S(0,z_2,\cdots,z_k,y_1,\cdots,y_k)|}f^k(x_1)\mathrm{d}x_1\mathrm{d}z_2\cdots\mathrm{d}z_k\mathrm{d}y_1\cdots\mathrm{d}y_k\nonumber
\end{multline}
by using the mean value theorem for integrals and the dominated
convergence theorem. Thereby, the first term on the right hand side
of (\ref{2}) tends to $k^{-1}\int_{\mathbb{R}^2}\varphi_T(\lambda
f(x_1))f(x_1)\mathrm{d}x_1$.

On the other hand, note that $n^{-1}{n \choose k}\le Cr_n^{-2(k-1)}$
for some positive constant $C$, then the absolute value of the
second term on the right hand side of (\ref{2}) is bounded by
$C\int_{\mathbb{R}^2}f(x_1)\eta_n(x_1)\mathrm{d}x_1$, where
$\eta_n(x_1)$ is defined by
$\eta_n(x_1)=\int_{D(x_1,kr_n)}\cdots\int_{D(x_1,kr_n)}r_n^{-2(k-1)}\\
\cdot\big|\prod_{i=2}^kf(x_i)-f^{k-1}(x_1)\big|\mathrm{d}x_2\cdots\mathrm{d}x_k$.
By the mean value theorem for integrals, $\eta_n(x_1)$ tends to $0$
if $x_1$ is a continuous point of $f$, then by the dominated
convergence theorem, as in the proof of Proposition 1,
$\int_{\mathbb{R}^2}f(x_1)\eta_n(x_1)\mathrm{d}x_1\rightarrow0$ as
$n\rightarrow\infty$. Thus the proof is completed. $\Box$

\medskip
\noindent\textbf{Proposition 4.}\quad\itshape For $k=1$, let $T$ be
a graph of order $k$, that is, $T$ is a single point. Suppose
$nR_n^d\rightarrow\lambda\in[0,\infty)$ in probability, as
$n\rightarrow\infty$. Then
$$
\lim_{n\rightarrow\infty}n^{-1}E\tilde{G}_{n,T}=\int_{\mathbb{R}^d}e^{-f(x)\theta\lambda}f(x)\mathrm{d}x
$$
\normalfont

\medskip
\noindent\textbf{Proof}. By the definition of $g_{n,T}$,
\begin{eqnarray*}
n^{-1}E\tilde{G}_{n,T}&=&\int_0^{\infty}\int_{\mathbb{R}^d}g_{n,T}(x_1,r_{n,1})\Big(1-\int_{B(x_1,r_{n,1})}f(x)\mathrm{d}x\Big)^{n-1}f(x_1)q_n(r_{n,1})\mathrm{d}x_1\mathrm{d}r_{n,1}\\
 &=&\int_0^{\infty}\int_{\mathbb{R}^d}\Big(1-\int_{B(x_1,r_{n,1})}f(x)\mathrm{d}x\Big)^{n-1}f(x_1)q_n(r_{n,1})\mathrm{d}x_1\mathrm{d}r_{n,1}\\
 &\sim&\int_0^{\infty}\int_{\mathbb{R}^d}e^{-nf(x_1)\theta
 r_{n,1}^d}f(x_1)q_n(r_{n,1})\mathrm{d}x_1\mathrm{d}r_{n,1}.
\end{eqnarray*}
Hence using the dominated convergence theorem and the assumption of
$R_n$, the above expression tends to
$\int_{\mathbb{R}^d}e^{-f(x_1)\theta\lambda}f(x_1)\mathrm{d}x_1$ as
$n$ tends to infinity, which concludes the proof. $\Box$

Now we are in position to prove our strong laws of large numbers.

\medskip
\noindent\textbf{Proof of Theorem 1}. In order to use Lemma 1, we
shall first define a filtration. Let
$\mathcal{F}_0=\{\emptyset,\Omega\}$ be the trivial $\sigma$-field,
and $\mathcal{F}_i=\sigma\{(X_1,Y_1),\cdots,(X_i,Y_i)\}$ for $i\ge
1$. Define a martingale difference sequence as
$D_{n,i}:=E(\tilde{H}_{n,T}|\mathcal{F}_i)-E(\tilde{H}_{n,T}|\mathcal{F}_{i-1})$,
therefore we can write
$\tilde{H}_{n,T}-E\tilde{H}_{n,T}=\sum_{i=1}^nD_{n,i}$. Let
$\tilde{H}_{n,T}^i$ denote the number of isolated $T$-subgraphs in
$G(\mathcal{X}_{n+1}\backslash\{X_i\},\mathcal{Y}_{n+1}\backslash\{Y_i\},r_n)$
and $\tilde{H}'_{n+1,T}$ denote the number of isolated $T$-subgraphs
in $G(\mathcal{X}_{n+1},\mathcal{Y}_{n+1},r_n)$. Thus we have
$D_{n,i}=E(\tilde{H}_{n,T}-\tilde{H}_{n,T}^i|\mathcal{F}_i)$. Adding
a point to a finite set in $\mathbb{R}^d$ can cause the number of
isolated $T$-subgraphs to increase by at most 1, and can cause it to
decrease by less than a constant $M$ (here $d=2$, so we may take
$M=6$), thereby we get
\begin{eqnarray*}
|\tilde{H}_{n,T}-\tilde{H}_{n,T}^i|&\le&|\tilde{H}_{n,T}-\tilde{H}'_{n+1,T}|+|\tilde{H}'_{n+1,T}-\tilde{H}_{n,T}^i|\\
 &\le&(M+1)+(M+1)=2(M+1)
\end{eqnarray*}
and then $D_{n,i}\le 2(M+1)$. Now for any $\varepsilon>0$, by Lemma
1, we have $P(|\tilde{H}_{n,T}-E\tilde{H}_{n,T}|>\varepsilon n)\le 2
\exp\big(-\frac{\varepsilon^2 n}{8(M+1)^2}\big)$, which is summable
in $n$. The result follows by Borel-Cantelli Lemma and Proposition
3.    $\Box$

\medskip
\noindent\textbf{Proof of Theorem 2}. The proof parallels to that of
Theorem 1. Define a filtration:
$\mathcal{F}_{n,0}=\{\emptyset,\Omega\}$ and
$\mathcal{F}_{n,i}=\sigma\{(X_1,R_{n,1}),\cdots,(X_i,R_{n,i})\}$ for
$i\ge 1$. A martingale difference sequence is defined by
$D_{n,i}=E(\tilde{G}_{n,T}|\mathcal{F}_{n,i})-E(\tilde{G}_{n,T}|\mathcal{F}_{n,i-1})$,
and then we have
$\tilde{G}_{n,T}-E\tilde{G}_{n,T}=\sum_{i=1}^nD_{n,i}$. Let
$\tilde{G}_{n,T}^i$ denote the number of isolated vertices in
$G(\mathcal{X}_{n+1}\backslash\{X_i\},\mathcal{R}_{n}\cup R_{n,n+1}
\backslash\{R_{n,i}\})$. It's easy to see that
$D_{n,i}=E(\tilde{G}_{n,T}-\tilde{G}_{n,T}^i|\mathcal{F}_{n,i})$.
Reason similarly as in the proof of Theorem 1, there exists a
constant $M>0$ depends only on $d$ such that
$|\tilde{G}_{n,T}-\tilde{G}_{n,T}^i|\le M$, hence $|D_{n,i}|\le M$
a.s. For $\varepsilon>0$, by Lemma 1,
$P(|\tilde{G}_{n,T}-E\tilde{G}_{n,T}|>\varepsilon n)\le 2
\exp\big(-\frac{\varepsilon^2 n}{2M^2}\big)$, which is summable in
$n$. The result then follows by Borel-Cantelli Lemma and Proposition
4. $\Box$

\medskip
\noindent\textbf{Proof of Theorem 3}. Proceeding along the same line
as the proof of Theorem 1, by Lemma 1, we get
$$
P\big(|\tilde{H}_{n,T}-E\tilde{H}_{n,T}|>\varepsilon
n^kr_n^{2(k-1)}\big)\le 2 \exp\bigg(-\frac{\varepsilon^2
n^{2k-1}r_n^{4(k-1)}}{8(M+1)^2}\bigg),
$$
which is summable in $n$ by the assumption. The result follows by
Borel-Cantelli Lemma and Proposition 2. $\Box$

\medskip
In the next proof, we substitute Lemma 2 for Lemma 1 since this time
the basic Azuma' inequality is no longer valid.

\medskip
\noindent\textbf{Proof of Theorem 4}. Set
$\eta:=\frac{1\wedge\delta}{3(k-1)}$ and partition $\mathbb{R}^2$
into squares ($A_{n,i}$, $i\in\mathbb{N}$), each of side $r_n$. Let
$$
E:=\{\mathcal{X}\subset\mathbb{R}^2|\mathrm{card}(\mathcal{X})=n,
\mathcal{X}(A_{n,i})\le n^{\eta}(nr_n^2\vee1)\ \mathrm{for}\
\mathrm{every}\ A_{n,i}\ s.t.\ A_{n,i}\cap
\mathrm{supp}f\not=\emptyset\}.
$$
Since $\mathcal{X}_n(A_{n,i})\sim
Bin(n,\int_{A_{n,i}}f(x)\mathrm{d}x)$ is a binomial random variable,
by a Chernoff bound (see e.g.\cite{2} p.16), we have
$P(\mathcal{X}_n(A_{n,i})>n^{\eta}(nr_n^2\vee1))\le e^{-n^{\eta}}$
when $n$ is large enough. Since $\mathrm{supp}f$ is bounded, we have
\begin{equation}
P(\mathcal{X}_n\not\in E)\le
Cn^{\frac{2k-1-\delta}{2k-2}}e^{-n^{\eta}}\label{3},
\end{equation}
where C is some positive constant. Now let's define a filtration.
Let $\mathcal{F}_0=\{\emptyset,\Omega\}$ be the trivial
$\sigma$-field, and
$\mathcal{F}_i=\sigma\{(X_1,Y_1),\cdots,(X_i,Y_i)\}$ for $i\ge 1$.
We have $H_{n,T}-EH_{n,T}=\sum_{i=1}^nD_{n,i}$ with the martingale
differences
$D_{n,i}:=E(H_{n,T}|\mathcal{F}_i)-E(H_{n,T}|\mathcal{F}_{i-1})$.
Let $H_{n,T}^i$ denote the number of induced $T$-subgraphs in
$G(\mathcal{X}_{n+1}\backslash\{X_i\},\mathcal{Y}_{n+1}\backslash\{Y_i\},r_n)$,
then $D_{n,i}=E(H_{n,T}-H_{n,T}^i|\mathcal{F}_i)$. Set an event
$E_{n,i}=\{\mathcal{X}_n\in E,\mathcal{X}_{n+1}\backslash\{X_i\}\in
E\}$, hence we may derive $|H_{n,T}-H_{n,T}^i|\cdot1_{E_{n,i}}\le
C\cdot[n^{\eta}(nr_n^2\vee1)]^{k-1}$ for some constant $C$, and
$|H_{n,T}-H_{n,T}^i|\le n^k$, since changing the position of one
point in a configuration can only at most affect the subgraphs
constructed by points in surrounding nine squares. Consequently,
\begin{eqnarray}
|D_{n,i}|&\le&E(|H_{n,T}-H_{n,T}^i|\cdot1_{E_{n,i}}|\mathcal{F}_i)+E(|H_{n,T}-H_{n,T}^i|\cdot1_{E_{n,i}^c}|\mathcal{F}_i)\nonumber\\
 &\le&C\cdot[n^{\eta}(nr_n^2\vee1)]^{k-1}+n^kP(E_{n,i}^c|\mathcal{F}_i)\label{4}
\end{eqnarray}
Define an event $F_{n,i}:=\{P(E_{n,i}^c|\mathcal{F}_i)\le n^{-k}\}$,
thus
$$
P(F_{n,i}^c)\le n^kE(P(E_{n,i}^c|\mathcal{F}_i))=n^kP(E_{n,i}^c)\le
C'n^{k+\frac{2k-1-\delta}{2k-2}}e^{-n^{\eta}},
$$
by Markov's inequality and (\ref{3}). Note by (\ref{4}) when
$F_{n,i}$ occurs, $|D_{n,i}|\le
C\cdot[n^{\eta}(nr_n^2\vee1)]^{k-1}$. Therefore, by using Lemma 2,
we obtain for any $\varepsilon>0$
\begin{eqnarray*}
P(|H_{n,T}-EH_{n,T}|>\varepsilon
n^kr_n^{2(k-1)})&\hspace{-5pt}\le&\hspace{-5pt}2\exp\bigg(-\frac{\varepsilon^2n^{2k}r_n^{4(k-1)}}{Cn[n^{\eta}(nr_n^2\vee1)]^{2(k-1)}}\bigg)\\
 & &+\bigg(1+\frac{2n^k}{\varepsilon
n^kr_n^{2(k-1)}}\bigg)\cdot C'n^{k+1+\frac{2k-1-\delta}{2k-2}}e^{-n^{\eta}}\\
 &\hspace{-5pt}\le&\hspace{-5pt}2e^{-\frac{\varepsilon^2}{C}\big(n^{1-2\eta(k-1)}\wedge
 n^{2k-1-2\eta(k-1)}r_n^{4(k-1)}\big)}+C''n^{3k-\delta}e^{-n^{\eta}},
\end{eqnarray*}
which is summable in $n$. The result then follows by Borel-Cantelli
Lemma and Proposition 1. $\Box$

\bigskip
\noindent{\Large\textbf{4. Conclusions and future work}}

\smallskip
We have derived various strong laws of large numbers of subgraph
counts in two types of directed random geometric networks. These
models of directed graphs are natural generalizations of standard
geometric graphs applicable to analysis of a variety of spatial
networks exemplified as wireless communication networks and sensor
networks. Our results obtained for the second model is limited and a
very natural question would be to ask what happens for other sorts
of small subgraphs.

\end{document}